\documentclass[10pt]{amsart}
\usepackage{amsmath,amsthm,amssymb}
\usepackage{graphicx}
\usepackage[margin=1.4in]{geometry}

\begin{document}


\newcommand{\PGL} {\Pj\Gl_2(\R)}                           
\newcommand{\PGLC} {\Pj\Gl_2(\C)}                          
\newcommand{\RP} {\R\Pj^1}                                 
\newcommand{\CP} {\C\Pj^1}                                 
\newcommand{\C} {{\mathbb C}}                              
\newcommand{\R} {{\mathbb R}}                              
\newcommand{\Z} {{\mathbb Z}}                              
\newcommand{\Pj} {{\mathbb P}}                             
\newcommand{\Sg} {{\mathbb S}}                             
\newcommand{\Gl} {{\rm Gl}}                                

\newcommand{\suchthat} {\:\: | \:\:}
\newcommand{\ore} {\ \ {\it or} \ \ }
\newcommand{\oand} {\ \ {\it and} \ \ }

\newcommand{\oM} [1] {\ensuremath{{\mathcal M}_{0,#1}(\R)}}                 
\newcommand{\M} [1] {\ensuremath{{\overline{\mathcal M}}{_{0, #1}(\R)}}}    
\newcommand{\cM} [1] {\ensuremath{{\mathcal M}_{0,#1}}}                     
\newcommand{\CM} [1] {\ensuremath{{\overline{\mathcal M}}{_{0, #1}}}}       

\newcommand{\Tubeset} {\mathfrak{T}}                                        

\newcommand{\D} {\Delta}
\newcommand{\K} {\mathcal{K}}                           
\newcommand{\J} {\mathcal{J}}

\newcommand{\sm}{\varepsilon}

\newcommand{\rec}[1] {G^*(#1)}                  

\newcommand{\KG} {{\K} G}
\newcommand{\JG} {{\J} G}                

\newcommand{\temp} {\nabla}
\newcommand{\JGr} {\JG_r}
\newcommand{\JGd} {\JG_d}

\newcommand{\tubeA}{\raisebox{-.5ex}{\includegraphics{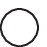}}}
\newcommand{\tubeB}{\raisebox{-.6ex}{\includegraphics{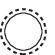}}}
\newcommand{\tubeC}{\raisebox{-.6ex}{\includegraphics{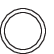}}}

%
%

\theoremstyle{plain}
\newtheorem{thm}{Theorem}
\newtheorem{prop}[thm]{Proposition}
\newtheorem{cor}[thm]{Corollary}
\newtheorem{lem}[thm]{Lemma}
\newtheorem{conj}[thm]{Conjecture}

\theoremstyle{definition}
\newtheorem{defn}[thm]{Definition}
\newtheorem*{exmp}{Example}

\theoremstyle{remark}
\newtheorem*{rem}{Remark}
\newtheorem*{hnote}{Historical Note}
\newtheorem*{nota}{Notation}
\newtheorem*{ack}{Acknowledgments}
\numberwithin{equation}{section}

\title {Marked tubes and the graph multiplihedron}


\author{Satyan L.\ Devadoss}
\address{S.\ Devadoss: Williams College, Williamstown, MA 01267}
\email{satyan.devadoss@williams.edu}

\author{Stefan Forcey}
\address{S.\ Forcey: Tennessee State University, Nashville, TN 37209}
\email{sforcey@tnstate.edu}

\begin{abstract}
Given a graph $G$, we construct a convex polytope whose face poset is based on marked subgraphs of $G$.  Dubbed the \emph{graph multiplihedron}, we provide a realization using integer coordinates.  Not only does this yield a natural generalization of the multiphihedron, but features of this polytope appear in works related to quilted disks, bordered Riemann surfaces, and operadic structures.  Certain examples of graph multiplihedra are related to Minkowski sums of simplices and cubes and others to the permutohedron.
\end{abstract}


\maketitle


\baselineskip=15pt

%
%
\section{Introduction}
\subsection{}

\begin{figure}[h]
\includegraphics[width=\textwidth]{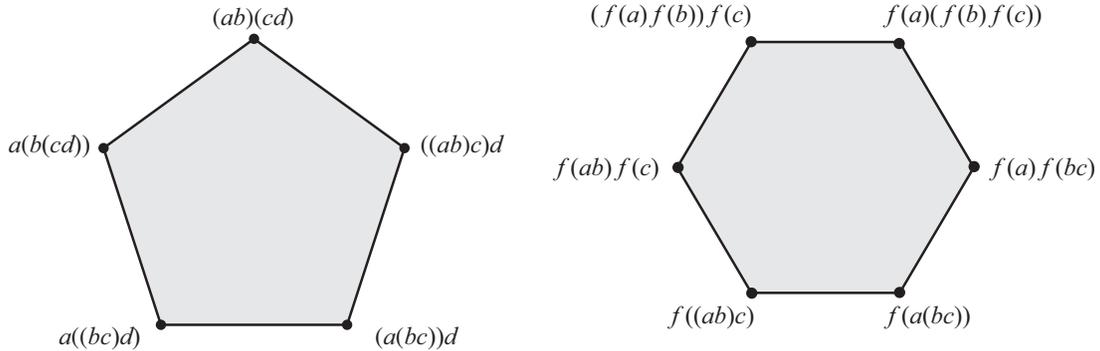}
\caption{The two-dimensional (a) associahedron $K_4$ and (b) multiplihedron $\J_3$.}
\label{f:2dmulti}
\end{figure}

The associahedron has continued to appear in a vast number of mathematical fields since its debut in homotopy theory \cite{sta}.
Stasheff classically defined the associahedron $K_n$ as a CW-ball with codim $k$ faces corresponding to using $k$ sets of parentheses meaningfully on $n$ letters; Figure~\ref{f:2dmulti}(a) shows the picture of $K_4$.
Indeed, the associahedron appears as a tile of \M{n}, the compactification of the real moduli space of punctured Riemann spheres \cite{dev1}.   Given a graph $G$, the \emph{graph associahedron} $\KG$ is a convex polytope generalizing the associahedron, with a face poset based on the connected subgraphs of $G$ \cite{cd}.  For instance, when $G$ is a path, a cycle, or a complete graph, $\KG$ results in the associahedron, cyclohedron, and permutohedron, respectively.  In \cite{dev2}, a geometric realization of $\KG$ is given, constructing this polytope from truncations of the simplex.  Indeed,  $\KG$ appears as tilings of minimal blow-ups of certain Coxeter complexes \cite{cd}, which themselves are natural generalizations of the moduli spaces \M{n}.

Our interests in this paper lie with the \emph{multiplihedron} $\J_n$, a polytope introduced by Stasheff in order to define $A_\infty$ maps between $A_\infty$ spaces \cite{sta2}.  Boardman and Vogt \cite{bv} fleshed out the definition in terms of \emph{painted trees};  a detailed combinatorial description was then given by Iwase and Mimura \cite{im}.  Saneblidze and Umble relate the multiplihedron to co-bar constructions of category theory and the notion of permutahedral sets \cite{su}.  In particular, $\J_n$ is a polytope of dimension $n-1$ whose vertices correspond to the ways of bracketing $n$ variables and applying a morphism $f$ (seen as an $A_\infty$ map).  Figure~\ref{f:2dmulti}(b) shows the two-dimensional hexagon which is $\J_3$.   Recently, Forcey \cite{for2} has provided a realization of the multiplihedron, establishing it as a convex polytope.  Moreover, Mau and Woodward \cite{mw} have shown $\J_n$ as the compactification of the moduli space of \emph{quilted disks}.

\subsection{}

In this paper, we generalize the multiplihedron to \emph{graph} multiplihedra $\JG$.  Indeed, the graph multiplihedra are already beginning to appear in literature; for instance, in \cite{dhv}, they arise as realizations of certain bordered Riemann disks of Liu \cite{liu}.  Similar to multiplihedra,  the graph multiplihedra degenerates into two natural polytopes;  these polytopes are akin to one measuring associativity in the domain of the morphism $f$ and the other in the range \cite{for1}.

An overview of the paper is as follows:  Section~\ref{s:defns} describes the graph multiplihedron as a convex polytope based on marked tubes, given by Theorem~\ref{t:poset}.  Section~\ref{s:geom} then follows with numerous examples.  When $G$ is a graph with no edges, we relate $\JG$ to Minkowski sums of cubes and simplices; when $G$ is a complete graph, $\JG$ appears as the permutohedron, the only graph multiplihedron which is a \emph{simple} polytope.  In Section~\ref{s:facets}, geometric properties of the facets of graph multiplihedra are discussed.  A realization of $\JG$ with integer coordinates is introduced in Section~\ref{s:real} along with constructions of two related polytopes.  Finally, the proof of the key theorems are provided in Section~\ref{s:proof}.

%
%
\section{Definitions} \label{s:defns}
\subsection{}

We begin with motivating definitions of graph associahedra; the reader is encouraged to see \cite[Section 1]{cd} for details.

\begin{defn}
Let $G$ be a finite graph.  A \emph{tube} is a set of nodes of $G$ whose induced graph is a connected subgraph of $G$.  Two tubes $u_1$ and $u_2$ may interact on the graph as follows:
\begin{enumerate}
\item Tubes are \emph{nested} if  $u_1 \subset u_2$.
\item Tubes \emph{intersect} if $u_1 \cap u_2 \neq \emptyset$ and $u_1 \not\subset u_2$ and $u_2 \not\subset u_1$.
\item Tubes are \emph{adjacent} if $u_1 \cap u_2 = \emptyset$ and $u_1 \cup u_2$ is a tube in $G$.
\end{enumerate}
Tubes are \emph{compatible} if they do not intersect and they are not adjacent.  A \emph{tubing} $U$ of $G$ is a set of tubes of $G$ such that every pair of tubes in $U$ is compatible.
\end{defn}

\begin{rem}
For the sake of clarity, a slight alteration of this definition is needed.  Henceforth, the entire graph (whether it be connected or not) will itself be considered a tube, called the \emph{universal tube}.  Thus all other tubes of $G$ will be nested within this tube.  Moreover, we force every tubing of $G$ to contain (by default) its universal tube.
\end{rem}

\begin{figure}[h]
\includegraphics[width=\textwidth]{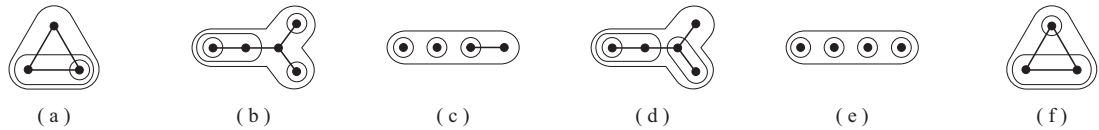}
\caption{(a)-(c) Allowable tubings and (d)-(f) forbidden tubings.}
\label{f:legaltubing}
\end{figure}

When $G$ is a disconnected graph with connected components $G_1$, \ldots, $G_k$, an additional condition is needed: If $u_i$ is the tube of $G$ whose induced graph is $G_i$, then any tubing of $G$ cannot contain all of the tubes $\{u_1, \ldots, u_k\}$.  Thus, for a graph $G$ with $n$ nodes, a tubing of $G$ can at most contain $n$ tubes.  Parts (a)-(c) of Figure~\ref{f:legaltubing} shows examples of allowable tubings, whereas (d)-(f) depict the forbidden ones.

\begin{thm} {\textup{\cite[Section 3]{cd}}} \label{t:graph}
For a graph $G$ with $n$ nodes, the \emph{graph associahedron} $\KG$ is a simple, convex polytope of dimension $n-1$ whose face poset is isomorphic to the set of tubings of $G$, ordered such that $U \prec U'$ if $U$ is obtained from $U'$ by adding tubes.
\label{d:pg}
\end{thm}

\begin{figure}[h]
\includegraphics[width=\textwidth]{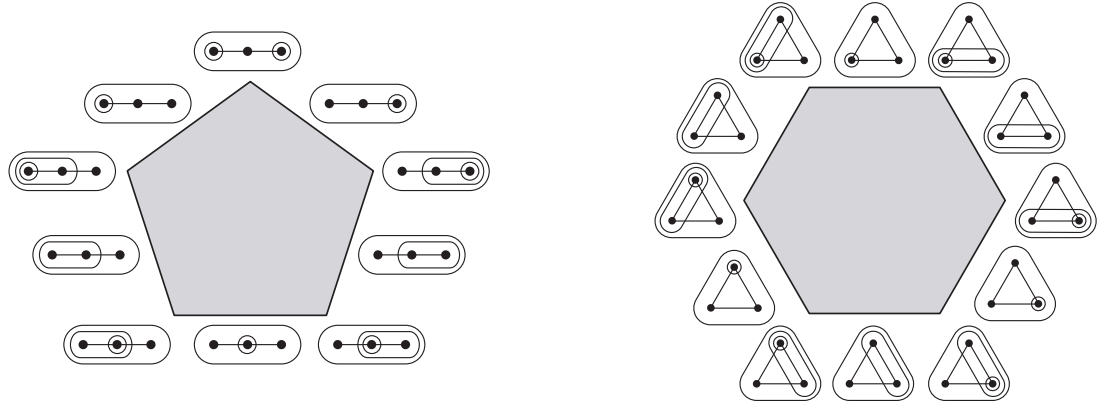}
\caption{Graph associahedra of a path and a cycle as underlying graphs.}
\label{f:kwexmp}
\end{figure}

\begin{exmp}
Figure~\ref{f:kwexmp} shows two examples of graph associahedra, having underlying graphs as paths and cycles, respectively, with three nodes.  These turn out to be the associahedron \cite{sta} and cyclohedron \cite{bt} polytopes.
\end{exmp}

\subsection{}
The notion of a tube is now extended to include markings.

\begin{defn}
A \emph{marked tube} of $G$ is a tube with one of three possible markings:
\begin{enumerate}
\item a \emph{thin} tube \ \tubeA \ \ given by a solid line,
\item a \emph{thick} tube \ \tubeC \ \ given by a double line, and
\item a \emph{broken} tube  \ \tubeB \ \ given by fragmented pieces.
\end{enumerate}
Marked tubes $u$ and $v$ are \emph{compatible} if
\begin{enumerate}
\item $u$ and $v$ do not intersect,
\item $u$ and $v$ are not adjacent, and
\item if $u \subset v$ where $v$ is not thick, then $u$ must be thin.
\end{enumerate}
A \emph{marked tubing} of $G$ is a collection of pairwise compatible marked tubes of $G$.
\end{defn}

\begin{figure}[h]
\includegraphics{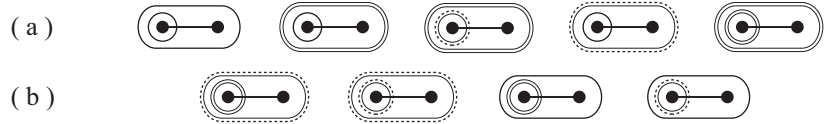}
\caption{(a) Allowable marked tubings and (b) forbidden marked tubings.}
\label{f:markedtubing}
\end{figure}

\noindent Figure~\ref{f:markedtubing} shows the nine possibilities of marking two nested tubes.  Out of these, row (a) shows allowable marked tubings, and row (b) shows those forbidden.

A partial order is now given on marked tubings of a graph $G$.  This poset structure is then used to construct the \emph{graph multiplihedron} below.  We start with a definition however.

\begin{defn}
Let $U$ be a tubing of graph $G$ containing tubes $u$ and $v$.  We say $u$ is \emph{closely nested} within $v$ if $u$ is nested within $v$ but not within any other tube of $U$ that is nested within $v$.  We denote this relationship as $u \Subset v$.
\end{defn}

\begin{figure}[h]
\begin{center}
\includegraphics[width=\textwidth]{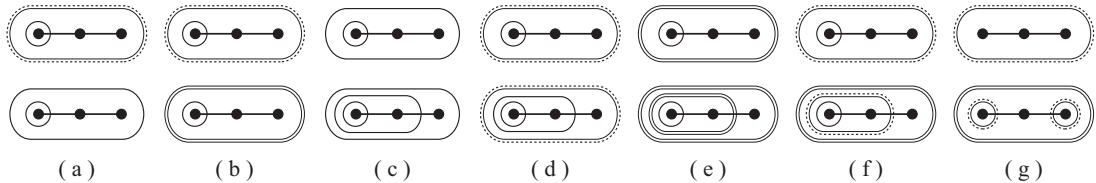}
\end{center}
\caption{The top row are the tubings and bottom row their refinements.}
\label{f:poset}
\end{figure}

\begin{defn} \label{d:poset}
The collection of marked tubings on a graph $G$ can be given the structure of a poset.   A marked tubings $U \prec U'$ if $U$ is obtained from $U'$ by a combination of the following three moves.  Figure~\ref{f:poset} provides the appropriate illustrations, with the top row depicting $U'$ and the bottom row $U$.
\begin{enumerate}
\item \emph{Resolving markings}:  A broken tube becomes either a thin tube (\ref{f:poset}a) or a thick tube (\ref{f:poset}b).
\item \emph{Adding thin tubes}:  A thin tube is added inside either a thin tube (\ref{f:poset}c) or broken tube (\ref{f:poset}d).
\item \emph{Adding thick tubes}:  A thick tube is added inside a thick tube (\ref{f:poset}e).
\item \emph{Adding broken tubes}:  A collection of compatible broken tubes $\{u_1, \ldots, u_n\}$ is added simultaneously inside a broken tube $v$ only when $u_i \Subset v$ \emph{and} $v$ becomes a thick tube; two examples are given in (\ref{f:poset}f) and (\ref{f:poset}g).
\end{enumerate}
\end{defn}

We are now in position to state one of our key theorems:

\begin{thm} \label{t:poset}
For a graph $G$ with $n$ nodes, the \emph{graph multiplihedron} $\JG$ is a convex polytope of dimension $n$ whose face poset is isomorphic to the set of marked tubings of $G$ with the poset structure given above.
\end{thm}

\begin{cor}
The codimension $k$ faces of $\JG$ correspond to marked tubings with exactly $k$ non-broken tubes.
\end{cor}

\noindent The proof of the theorem, along with the corollary, follows from the geometric realization of the graph multiplihedron given by Theorem~\ref{t:main}.  We postpone its proof until the end of the paper.

%
%
\section{Examples} \label{s:geom}
\subsection{}

\begin{figure}[h]
\includegraphics[width=.95\textwidth]{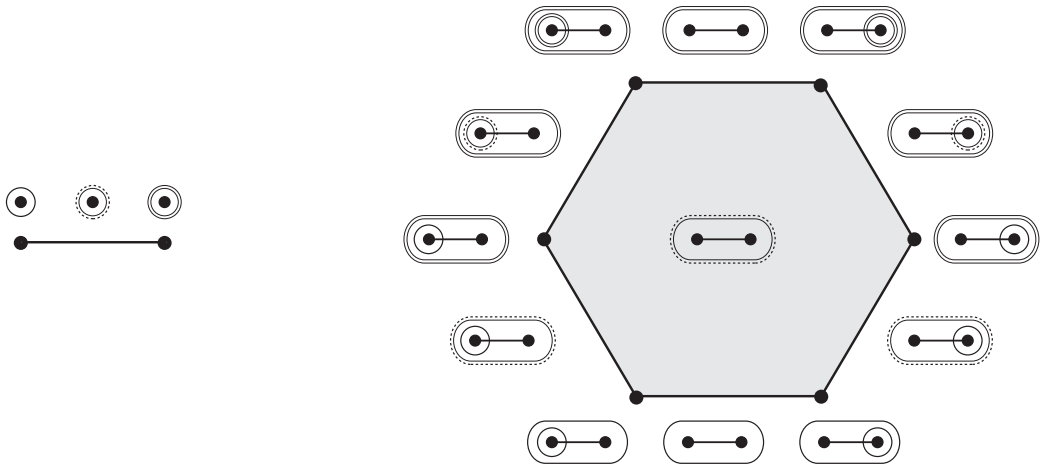}
\caption{The graph multiplihedron of a path with (a) one vertex $\J_2$ and (b) two vertices $\J_3$, along with labelings of faces by marked tubings.}
\label{f:path}
\end{figure}

The multiplihedron $\J_n$ serves as a parameter space for homotopy multiplicative morphisms.  From a certain perspective, as shown in \cite{su}, it naturally lies between the associahedron and the permutohedron.  If $G$ is a path with $n-1$ nodes, it is easy to see that $\JG$ produces the classical multiplihedron $\J_n$ of dimension $n-1$.  Figure~\ref{f:path}(a) shows the one-dimensional multiplihedron $\J_2$ as the interval, with endpoints labeled by legal tubings of a vertex.  The two-dimensional multiplihedron $\J_3$ is given in Figure~\ref{f:path}(b) with labeling by marked tubings; compare this with Figure~\ref{f:2dmulti}(b).  Notice that each vertex of $\JG$ corresponds to maximally resolved marked tubings, those with only thin or thick tubes.  The thick tubes capture multiplication in the domain of the morphism $f$, whereas the thin ones record the range.

\begin{figure}[h]
\includegraphics[width=\textwidth]{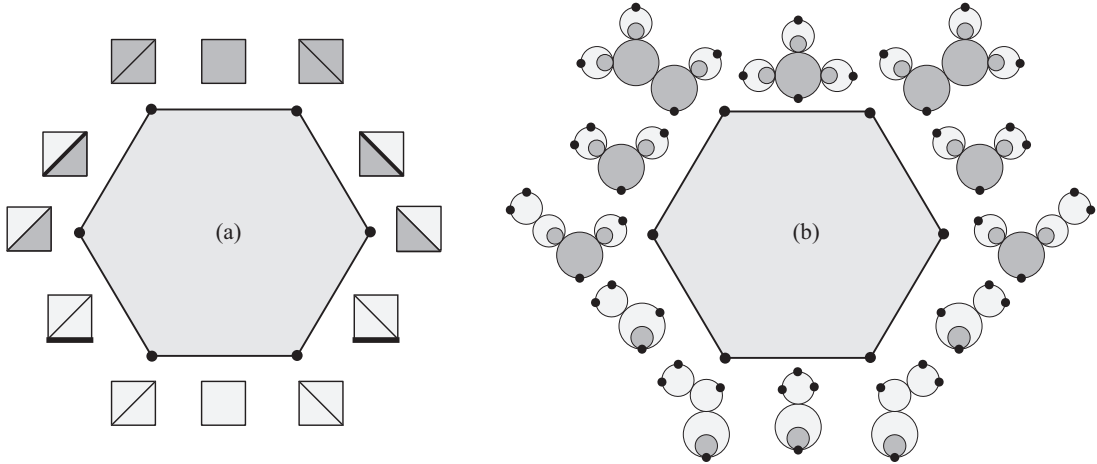}
\caption{The multiplihedron $J_3$ labelings by (a) painted diagonals of polygons and (b) quilted disks.}
\label{f:path-alt}
\end{figure}

Figure~\ref{f:path-alt} shows two different labelings of $\J_3$.  The left picture depicts the labeling using \emph{painted} diagonals of a polygon; these are dual to the painted trees of Boardman-Vogt \cite{bv} and Forcey \cite{for2}.  The right hexagon in Figure~\ref{f:path-alt} is labeled using the \emph{quilted disk} moduli spaces of Mau and Woodward \cite{mw}.  We leave it to the reader to construct bijections between these labelings of $\J_n$ and marked tubings on paths.

\subsection{}

\begin{figure}[h]
\includegraphics[width=\textwidth]{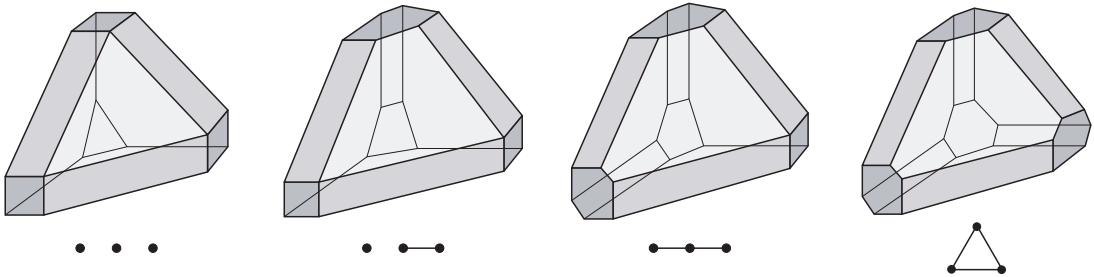}
\caption{The four possible three-dimensional graph multiplihedra.}
\label{f:mg-3d}
\end{figure}

There are only two kinds of graph multiplihedra $\JG$ when $G$ contains two nodes, one with $G$ disconnected and the other with $G$ being a path.  It is interesting to note that in both cases, $\JG$ is the hexagon, with labeling identical to Figure~\ref{f:path}.  This low-dimensional case is an anomaly, however.  Figure~\ref{f:mg-3d} shows the four different types of graph multiplihedra when $G$ contains three vertices. Notice that all of them but the rightmost polyhedron (when $G$ is a complete graph) are not simple.   Indeed, the rightmost graph multiplihedron of Figure~\ref{f:mg-3d} is combinatorially equivalent to the \emph{permutohedron}.  This is true in general as we now show.

\begin{thm} \label{t:perm}
Let $G$ be a complete graph on $n-1$ vertices.  The graph multiplihedron $\JG$ is combinatorially equivalent to the permutahedron $\mathcal{P}_n$.
\end{thm}

\begin{figure}[h]
\includegraphics{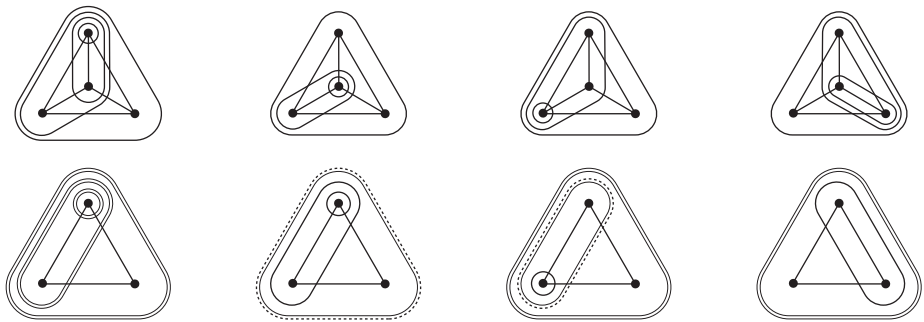}
\caption{Examples of tubings of $H$ (top row) and their images in $\psi_x$ (bottom row).}
\label{f:permutomap}
\end{figure}

\begin{proof}
Let $H$ be a complete graph on $n$ vertices and let $x$ be a node of $H$.
Let $G$ be the complete graph on $n-1$ vertices obtained from deleting $x$ from $H$.  We use the fact from \cite{dev2} that the permutohedron $\mathcal{P}_n$ is equivalent to the graph associahedron $\K H$.  Now we define a poset isomorphism $\psi_x$ from the $\K H$ (tubings of $H$) to $\JG$ (marked tubings of $G$).

Let $T$ be a tubing of $H$ (an element of $\K H$) and let $\sm(x)$ be the smallest tube of $T$ containing $x$.  Then $\psi_x(T)$ is the marked tubing of $G$ (an element of $\JG$), with tubes $\{u-x\}$ for tubes $u$ in $T$, where the marking of $u-x$ is as follows:
\begin{enumerate}
\item \emph{thick} if $\sm(x) \subset u$.
\item \emph{broken} if $\sm(x) = u$ and $u-x$ is not in $T$.
\item \emph{thin} otherwise.
\end{enumerate}
Figure~\ref{f:permutomap} shows four examples where the top row shows tubings of $H$ and the bottom row shows the image in $\psi_x$; in all four cases, $x$ is the top most point in the complete graph on four vertices.  Notice that if $\sm(x)=u$, the marked tubing $u-x$ in $\psi_x(T)$ will be broken only if there is another node whose smallest tube is $u$; otherwise $\sm(x)-x$ will be the same as some tube $u'$ not containing $x$. With these facts in mind, it is straightforward to check that $\psi_x$ is an isomorphism of posets.
\end{proof}

\begin{cor}
The graph multiplihedron is a \emph{simple} polytope only when $G$ is a complete graph.
\end{cor}

\begin{proof}
Let $G$ not be a complete graph, and let $a, b$ be two of the $n$ nodes of $G$ not connected by an edge.  Consider a maximal marked tubing $T$ on $G$ (corresponding to some vertex of $\JG$) consisting of $n$ thick tubes, two of which are the precisely the tubes $\{a\}$ and $\{b\}$.  We claim there are at least $n+1$ marked tubings $S$ such that $T \prec S$ and there exists no other tubing $S'$ where $T \prec S' \prec S$.  Find $n-1$ of them by removing any of the thick tubes except the universal one. Find the other two by making $\{a\}$ or $\{b\}$ into a broken tube. Thus the vertex labeled by $T$ in contained in  at least $n+1$ edges, so $\JG$ is not simple. The converse follows from the pervious theorem.
\end{proof}

\subsection{}

Let $\temp_n$ denote the graph multiplihedron for the graph with $n$ disjoint nodes.  The right side of Figure~\ref{f:path} shows $\temp_2$, Figure~\ref{f:mg-3d}(a) displays $\temp_3$, and the left side of Figure~\ref{f:mg-4d} provides the four-dimensional $\temp_4$ polytope.  We show an alternate construction of $\temp_n$ using Minkowski sums.

\begin{defn}
The \emph{Minkowski sum} of two point sets $A$ and $B$ in $\R^n$ is
$$A \oplus B \ := \ \{x + y \suchthat x \in A, y \in B\},$$
where $x+y$ is the vector sum of the two points.
\end{defn}

\begin{figure}[h]
\includegraphics[width=\textwidth]{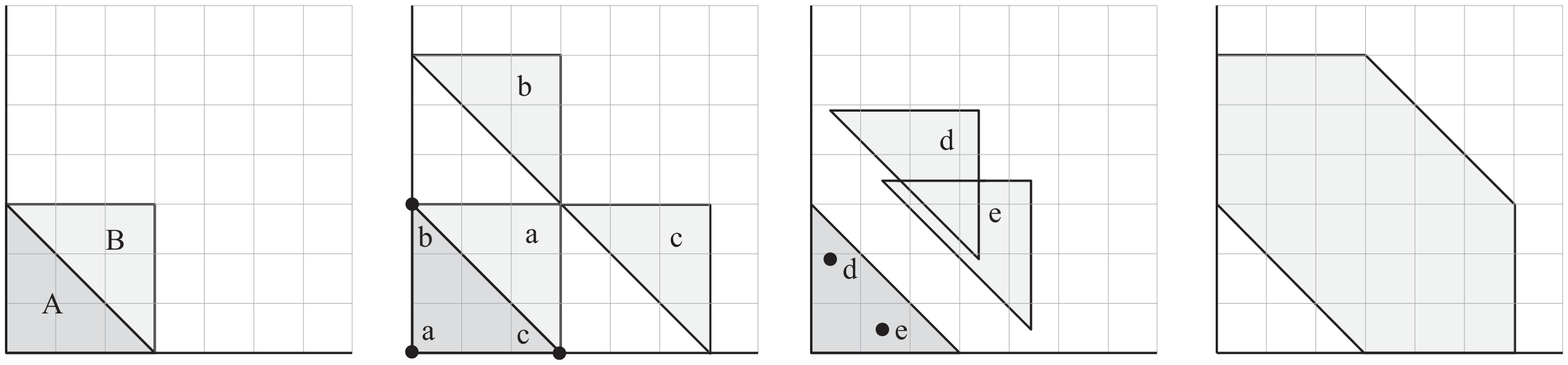}
\caption{The Minkowski sum of the two sets $A$ and $B$ on the left is given in the right.}
\label{f:minkowski}
\end{figure}

\begin{exmp}
The left side of Figure~\ref{f:minkowski} shows two sets $A$ and $B$, the decomposition of the square into two simplices.  The middle two figures display the sum of $B$ with certain labeled points of $A$, whereas the Minkowsi sum $A \oplus B$ is given in the right as the hexagon.
\end{exmp}

\begin{prop}
If $C_n$ is the $n$-cube $[0,1]^n$ in $\R^n$, then the hyperplane $\sum x_i = 1$ cuts $C_n$ into two polytopes, the simplex $\Delta_n$ and its complement $C_n - \Delta_n$.  The polytope $\temp_n$ is combinatorially equivalent to $\Delta_n \oplus (C_n - \Delta_n)$.
\end{prop}

\begin{proof}
We will demonstrate an isomorphism between vertex sets of the two polytopes (from $\temp_n$ to the above Minkowski sum) which preserves facet inclusion of vertices.  The vertices come in two groupings, and the bijection may be described piecewise on those sets.

{\textsc{Group I}}:  \ \ In the Minkowski sum,  the first grouping consists of $n$ vertices resulting from adding the origin to the vertices of the $(n-1)$-simplex facet of $C_n -\Delta_n.$  In $\temp_n$,  the first grouping consists of the $n$ vertices which correspond to the entirely thin maximal tubings; of course, for our edgeless graph, proper tubes consist of a single node.  For each node $v_i$, there is one of these tubings  which does not include $v_i$ itself as a tube. Map the vertex of $\temp_n$ corresponding to node $v_i$ not being a tube to the vertex of the Minkowski sum which lies on the $x_i$-axis.

{\textsc{Group II}}: \ \ The second grouping of vertices in $\temp_n$ consists of vertices with associated tubing containing the thick universal tube and all but one of the single nodes as either thick or thin tubes. Thus there are $n \cdot 2^{n-1}$ of these:  Choose the node that will not be a tube, then choose a (possibly empty) subset of the remaining nodes to be thick tubes.  The second grouping of vertices in the Minkowski sum are those resulting from a nonzero vertex of $\Delta_n$ being added to  a nonzero vertex of $C_n.$ The geometry dictates that for each facet of $C_n$ which is parallel to but not contained in a coordinate hyperplane, there will be $2^{n-1}$ vertices of the Minkowski sum --- one for each vertex of that facet of $C_n.$ These result from adding the vertex of $\Delta_n$ which lies in the axis perpendicular to the facet of $C_n$ to each of the vertices of that facet. Thus there are $n \cdot 2^{n-1}$ vertices in this second grouping. The bijection takes the vertex of $\temp_n$ associated to the tubing without the tube $v_p$ but with the thick tubes $v_{i_1}, \ldots, v_{i_k}$ to a vertex of the Minkowski sum formed by adding the vertex of $\Delta_n$ which lies in the $x_p$-axis to the vertex of $C_n$ which lies in the subspace spanned by the axes $\{x_p, x_{i_1}, \ldots, x_{i_k}\}$.

{\textsc{Facet Inclusion}}: \ \ We check that the bijection of vertices preserves facet inclusion. To check the first grouping of vertices, note that the $n+1$ lower facets of $\temp_n$ are given by a choice of a single thin single-node tube. The $n+1$ lower facets of the Minkowski sum correspond to adding the origin to a facet from $C_n - \Delta_n$, either to the facet which it shared with $\Delta_n$ or to one which lay in a coordinate hyperplane.  To check the second grouping of vertices, note first that in the Minkowski sum, the facets of $C_n - \Delta_n$ in the coordinate hyperplanes are extended by the vectors of $\Delta_n$ which lie in the same coordinate hyperplane.  Moreover, the $2^{n}-1$ upper facets of $\temp_n$ correspond to subsets of nodes  which will be the broken tubes.  The $2^{n}-1$ upper facets of the Minkowski sum correspond to adding the face of $C_n$ determined by intersecting a nonempty subset of the facets of $C_n$ which do not lay in a coordinate hyperplane to the orthogonal face of $\Delta_n$.  It is straightforward to verify that our bijection takes the vertices of a facet of $\temp_n$ to the vertices of a facet of the Minkowski sum.
\end{proof}

\begin{rem}
The construction of $\temp_2$ using this method is given in Figure~\ref{f:minkowski}.
\end{rem}

\begin{rem}
In \cite{pos} Postnikov defines the \emph{generalized permutohedra}, a class which encompasses a great many varieties of combinatorially defined polytopes. A subclass of these named \emph{nestohedra}, which include examples such as the graph associahedra and the Stanley-Pittman polytopes, are based on \emph{nested sets} as in Definition 7.3 of \cite{pos}. For the nestohedra, Postnikov gives a realization formulated as a Minkowski sum of simplices. A question deserving of further thought is whether there is a consistent definition of \emph{marked nested sets}, fitting into the scheme of generalized permutohedra, which would specialize to our marked tubings. It would be especially interesting to elucidate whether the Minkowski sum for $\temp_n$ discussed here has a nice generalization in that context.
\end{rem}

%
%
\section{Geometry of the Facets} \label{s:facets}
\subsection{}

The discussion and results in this section can be interpreted as describing either the poset of marked tubings of a graph $G$ or (after using Theorem~\ref{t:poset}) as describing the polytope which realizes this poset as its set of faces ordered by inclusion.  Therefore we will abuse notation and use $\JG$ to mean either the marked tubings themselves or the face poset labeled by them.   Our main concern here is regarding the \emph{facets} of $\JG$, the codimension one faces.  It follows immediately from the poset ordering given in Definition~\ref{d:poset} that the facets of $\JG$ are the tubings which contain exactly one unbroken tube.  We refine the facets further:

\begin{defn}
The facets can be partitioned into two classes:  The \emph{upper tubings} contain exactly one thick universal tube and the \emph{lower tubings} contain exactly one thin tube.\footnote{We abuse terminology by calling them upper and lower facets as well.}
\end{defn}

\begin{figure}[h]
\includegraphics[width=\textwidth]{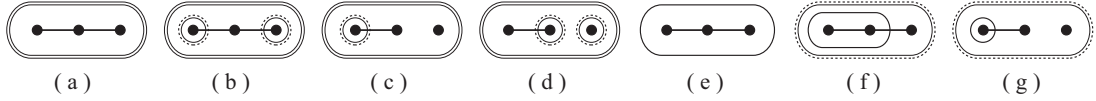}
\caption{Examples of (a-d) upper and (e-g) lower tubings.}
\label{f:upperlower}
\end{figure}

\noindent Figure~\ref{f:upperlower}(a-d) shows examples of upper tubings, whereas (e-g) show lower tubings.  Parts (a) and (e) show the universal thick and thin tubes, respectively.

\begin{lem} \label{l:enum}
Let $G$ be a graph with $n$ nodes.
\begin{enumerate}
\item The number of upper facets of $\JG$ is $2^{n-1}$.
\item The number of lower facets of $\JG$ equals one more than the number of facets of $\KG$.
\end{enumerate}
\end{lem}

\begin{proof}
The number of upper facets correspond to the number of ways to choose a nonempty set of nodes of $G$.   Since for each choice there is exactly one way to enclose the chosen nodes in a set of compatible broken tubes, we obtain $2^{n-1}$.  There exists a lower facet of $\JG$ (and a facet of $\KG$) for each tube of $G$.  However, $\JG$ has the additional lower facet corresponding to the thin universal tube. 
\end{proof}

\subsection{}

Before describing the geometry of the facets of $\JG$, a definition from \cite[Section 2]{cd} is needed. 

\begin{defn}
For graph $G$ and a collection of nodes $t$, construct a new graph $\rec{t}$ called the \emph{reconnected complement}: If $V$ is the set of nodes of $G$, then $V-t$ is the set of nodes of $\rec{t}$.  There is an edge between nodes $a$ and $b$ in $\rec{t}$ if either $\{a,b\}$ or $\{a,b\} \cup t$ is connected in $G$.
\end{defn}

\begin{figure}[h]
\includegraphics {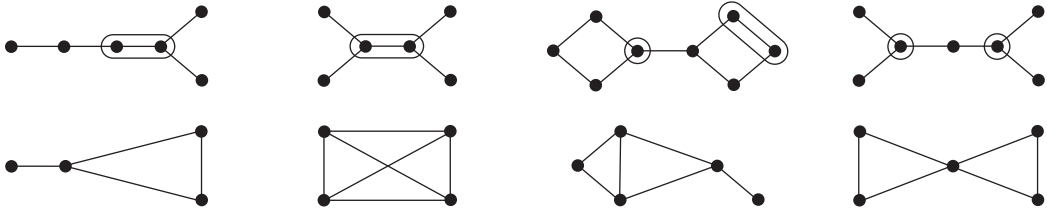}
\caption{Examples of tubes and their reconnected complements.}
\label{f:recon}
\end{figure}

\noindent Figure~\ref{f:recon} illustrates some examples on graphs along with their reconnected complements. For a given tube $t$ and a graph $G$, let $G(t)$ denote the induced subgraph on the graph $G$.  By abuse of notation, we sometimes refer to $G(t)$ as a tube.

\begin{prop} \label{p:lower}
Let $V$ be a lower facet of $\JG$ and let $t$ be the thin tube of $V$.  The face poset of $V$ is isomorphic to $\J \rec{t} \times \KG(t)$.  In particular, if $t$ is the universal thin tube, then $V$ is isomorphic to $\KG$.
\end{prop}

\begin{proof}
The last statement is easiest to verify. For the tubing $V$ consisting of the thin universal tube, any refinement of this tubing must be accomplished by adding more thin tubes. Thus the collection of refinements is just the poset of (all thin) tubings of $G$, and trivially isomorphic to $\KG$ by forgetting marking.

Now for the case in which the universal tube is broken, the marked tubings $U \prec V$ are those that contain the thin tube $t$.
First, any tubing $\{t_i\}$ of $\KG(t)$ becomes a marked tubing of $\JG$ by marking all the tubes as thin.  Let $\sigma(t_i)$ denote the marked tube of $\JG(t)$ achieved by assigning the thin marking.  Consider the map
$$\tau : \{\text{ marked tubes of } \rec{t}~\} \to \{\text{ marked tubes of $G$ containing $t$ }\}$$
given by
$$\tau(t') = \begin{cases}
    t'\cup t & \ \ \ \text{if  $t' \cup t$ is connected in $G$ or if $t' = \rec{t}$} \\
    t' & \ \ \ \text{otherwise}.
\end{cases}$$
Here $\tau(t')$ is defined to have the same marking as $t'$.  Now define a map
\begin{equation}
\hat{\tau}:\J \rec{t} \times \KG(t) \to V: (T,W) \mapsto \bigcup_{t_j \in T} \tau(t_j) \ \cup \ \bigcup_{t_i \in W}\sigma(t_i).
\label{e:tau}
\end{equation}
This is an isomorphism of posets by comparison to Theorem 2.9 of \cite{cd}.
\end{proof}

\begin{prop}\label{p:upper}
Let $V$ be an upper facet of $\JG$ and let $t_1,\dots,t_k$ be the broken tubes of $V$.  Let $t$ be the union of $\{t_i\}$.  The face poset of $V$ is isomorphic to
$$\K \rec{t} \times \JG(t_1) \times \cdots \times \JG(t_k).$$
In particular, if $V$ has no broken tubes, then $V$ is isomorphic to $\KG$.
\end{prop}

\begin{proof}
Again, we verify the last statement first. For the tubing $V$ with the only tube being the thick universal one, any refinement of this tubing must be accomplished by adding more thick tubes. Thus the collection of refinements is just the poset of (all thick) tubings of $G$, and isomorphic to $\KG$ by forgetting marking.

Let $t$ be the union of  the broken tubes $t_1,\dots,t_k$.  Consider the map
$$\eta : \{\text{ unmarked tubes of } \rec{t}~\} \to \{\text{ marked tubes of $G$ containing $t$ }\}$$
where
$$\eta(t') = \begin{cases}
     G & \ \ \ \text{ if $t'$ is universal } \\
     t' \ \cup \ \bigcup\{t_i ~|~ \exists~ u\in t' \text{ with } t_i \cup u \text{ connected }\} &  \ \ \ \text{otherwise}.
\end{cases}$$
Here $\eta(t')$ is defined to have the thick marking.
Now if $t''$ is a marked tube of $G(t_i)$ then $t''$ is also a marked tube of $G$.  Define a map
\begin{equation}
\hat{\eta}:\K \rec{t} \times \JG(t_1) \times \cdots \times \JG(t_k) \to V : (S,T_1,\dots,T_k) \mapsto \bigcup_{t' \in S} \eta(t') \ \cup \ \bigcup_{i= 1\dots k} T_i,
\label{e:eta}
\end{equation}
where the universal tube $G$ in the image is marked as thick.  This is poset isomorphism.
\end{proof}

\begin{rem}
The previous two Propositions can be seen as generalizations of the product structure of associahedra and cyclohedra as given in \cite{mss}.
\end{rem}

\begin{figure}[h]
\includegraphics[width=.95\textwidth]{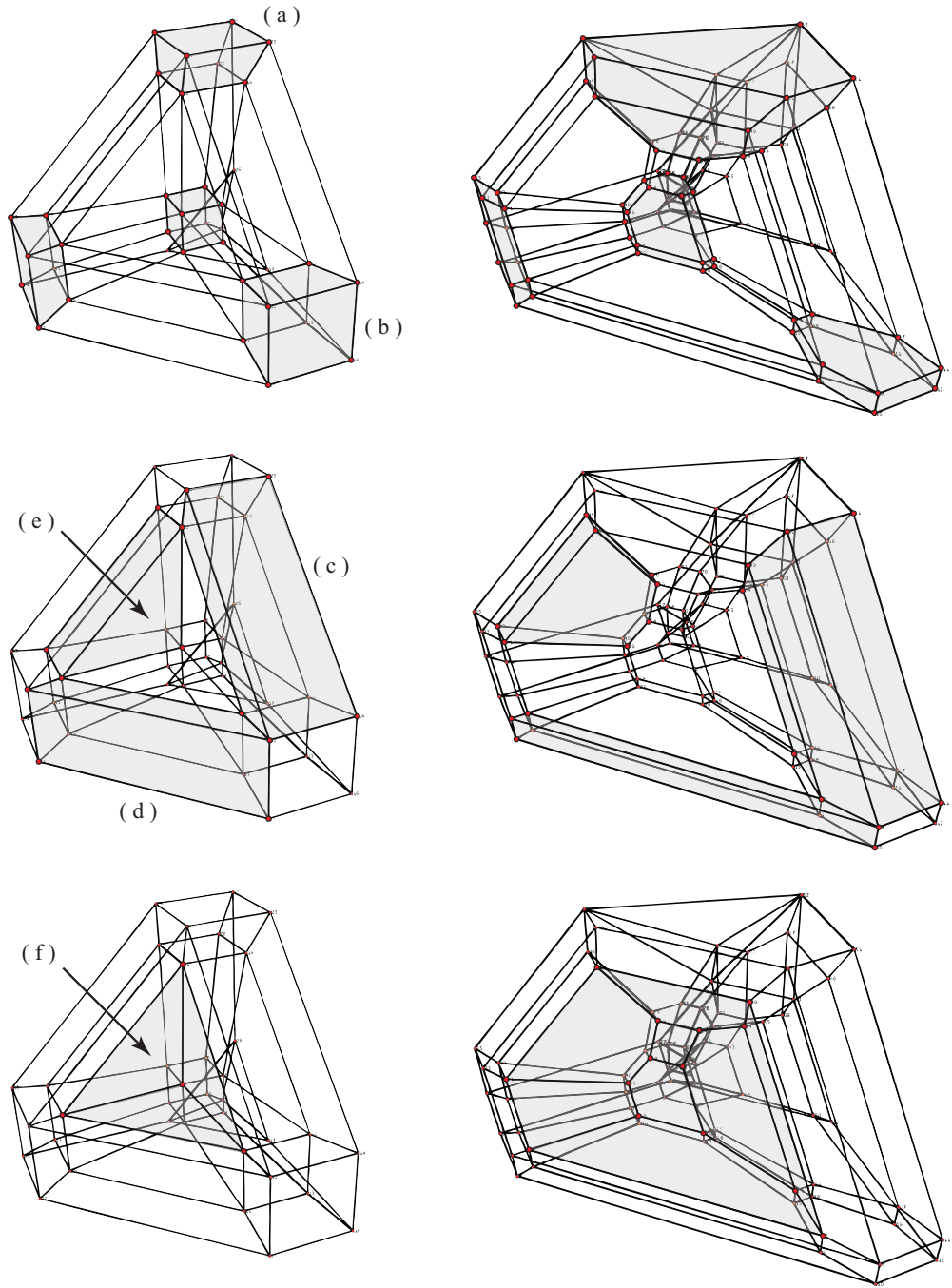}
\caption{The left side shows the Schlegel diagram of $\temp_4$ and the right side the classical four-dimensional multiplihedron.  Compare with Figure~\ref{f:mg-4d-pieces}.}
\label{f:mg-4d}
\end{figure}

\subsection{}

It turns out that the three-dimensional graph multiplihedra, all depicted in Figure~\ref{f:mg-3d}, do not yield complicated combinatorial structures.  It is not until four dimensions that certain ideas become transparent as given by Figure~\ref{f:mg-4d}.   The left side of this picture\footnote{The \textsc{Polymake} software \cite{pol} was used to construct these Schlegel diagrams with inputs of coordinates given by the realization of $\JG$ in Theorem~\ref{t:main}.} shows $\temp_4$ whereas the right side portrays the classical multiplihedron $\J_5$.   This perspective of the Schlegel diagram was chosen since the facets  visible are the upper facets.  Indeed, apparent from the two Propositions above, the complexity of $\JG$ is most prevalent in the structure of upper tubings.  Certain upper facets of $\temp_4$ are shaded here, with their corresponding facets similarly shaded on the right side.   

\begin{exmp}
Figure~\ref{f:mg-4d-pieces} analyzes the labeled facets in Figure~\ref{f:mg-4d}, the left side providing the geometry and the tubing label when $G$ has no edges and the right side when $G$ is a path.  The geometry of these upper facets of $\JG$ can be calculated using Proposition~\ref{p:upper}.  In what follows, understand that whenever the $n$-simplex $\Delta_n$ is mentioned, it arises from the graph \emph{associahedron} $\KG$, where $G$ is the graph with $n+1$ disjoint nodes.

\begin{figure}[h]
\includegraphics{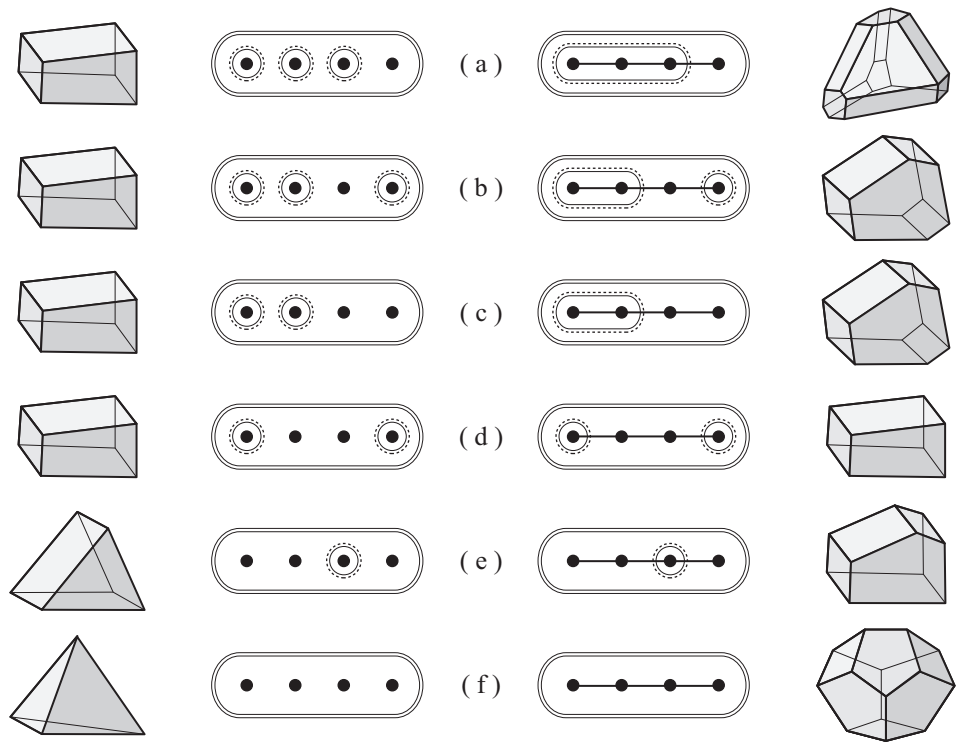}
\caption{Certain shaded upper facets of $\temp_4$ and $\J_5$ given in Figure~\ref{f:mg-4d} along with their labels by tubings.}
\label{f:mg-4d-pieces}
\end{figure}

\begin{enumerate}
\item[(a)]  On the left, the geometry of this upper facet is $\K_2 \times \J_2 \times \J_2 \times \J_2$.  Since $\K_2$ is a point and $\J_2$ is an edge, this is equivalent to a cube.  On the right, however, the three disjoint broken tubes combine into one large broken tube with three vertices.  Thus the geometry becomes $\K_2 \times \J_4$, resulting in the three-dimensional multiplihedron $\J_4$.
\item[(b)]  The left side labeling yields the same product structure as (a), resulting in a cube.  However, for the right side, one obtains $\K_2 \times \J_3 \times \J_2$, a hexagonal prism.
\item[(c)]  On the left, this facet is given by $\Delta_1 \times \J_2 \times \J_2$.  For the right side, we obtain $\K_3 \times \J_3$, a hexagonal prism like (b).  Note that although both (b) and (c) result in geometrically identical prisms, they encode different combinatorial data.
\item[(d)] Both kinds of facets are cubes.  The left is identical to part (c), whereas the right is $\K_3 \times \J_2 \times \J_2$.
\item[(e)] The left side labeling yields $\Delta_2 \times \J_2$, a triangular prism.  This transforms into a pentagonal prism $\K_4 \times \J_2$ on the right.
\item[(f)] The 3-simplex on the left becomes $\K_5$ on the right.
\end{enumerate}
\end{exmp}

%
%
\section{Realizations} \label{s:real}
\subsection{}

Thus far, our focus has been on the combinatorial structure of the graph multiplihedron based on marked tubings.  This section provides a geometric backbone giving $\JG$ a realization with integer coordinates.  Let $G$ be a graph with $n$ nodes, denoted $v_1, v_2, \ldots v_n$.  Let $M_G$  be the collection of maximal marked tubings of $G$.  Indeed, elements of $M_G$ will correspond to the vertices of $\JG$.  Notice that each tubing $U$ in $M_G$ contains exactly $n$ tubes, with each tube being either thin or thick. So $U$ assigns a unique tube $\sm(v)$ to each node $v$ of $G$, where $\sm(v)$ is the smallest tube in $U$ containing $v$.  Parts (a)-(c) of Figure~\ref{f:legaltubing} shows examples of maximal tubings of $G$.

\begin{figure}[h]
\includegraphics[width=\textwidth]{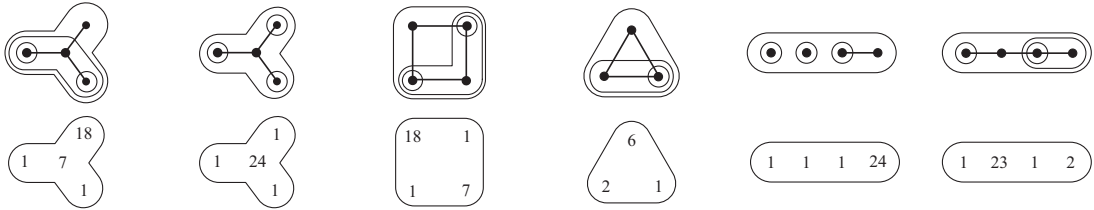}
\caption{Integer values of nodes associated to tubings.}
\label{f:coord}
\end{figure}

For each tubing $U$ in $M_G$, we define a function $f_U$ from the nodes of $G$ to the integers as follows:
$$f_U(v) = 3^{|\sm(v)|-1}-\sum_{s \Subset \sm(v) }3^{|s|-1}.$$
Note that $f_U$ is defined independently of the markings associated to the tubes of $U$.    Figure~\ref{f:coord} gives some examples of integer values of nodes associated to tubings.

\noindent Let $G$ be a graph with an ordering $v_1, v_2, \ldots, v_n$ of its nodes.  Define a map
\begin{equation}
c:M_G \to \R^n: U \mapsto (x_1,\dots,x_n)
\label{e:real}
\end{equation}
where
$$x_i =
\begin{cases}
\hspace{4.5mm} f_U(v_i) & \text{if $\sm(v_i)$ is thin.} \\
3 \cdot f_U(v_i)  & \text{if $\sm(v_i)$ is thick.}
\end{cases} $$

\noindent We are now in position to state the main theorem:

\begin{thm} \label{t:main}
For a graph $G$ with $n$ nodes, the convex hull of the points $c(M_G)$ in $\R^n$ yields the graph multiplihedron $\JG$.
\end{thm}

\begin{rem}
This theorem implies that the convex hull of the points $c(M_G)$ will produce a convex polytope whose face poset structure is given by $\JG$.  Thus, this geometric result immediately implies the combinatorial result of Theorem~\ref{t:poset}.  The proof of this theorem is given in Section~\ref{s:proof}, at the end of the paper.
\end{rem}

\begin{figure}[h]
\includegraphics[width=.9\textwidth]{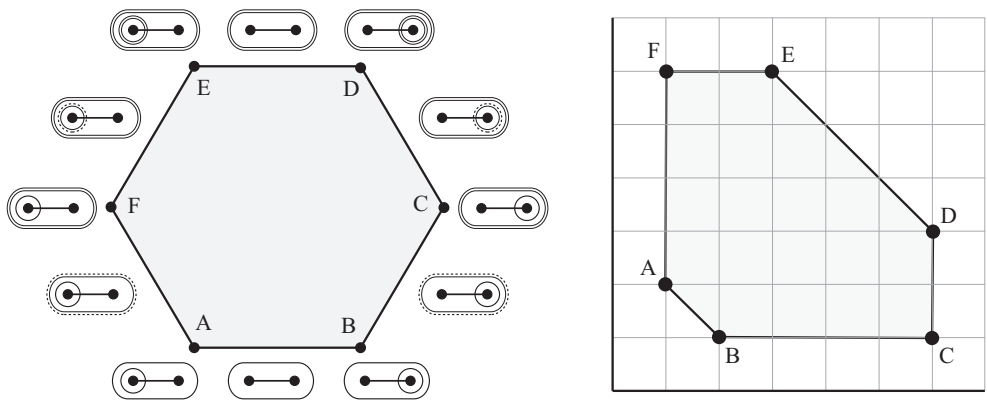}
\caption{The graph multiplihedron of a path along with its realization.}
\label{f:path-2d}
\end{figure}

\begin{exmp}
Figure~\ref{f:path-2d} shows an example of this realization.  The left side displays the hexagon poset given in Figure~\ref{f:path}, along with labels of the vertices.  The right side constructs the convex hull using integer coordinates on $\R^2$, with appropriate labelings of the vertices.
\end{exmp}

\subsection{}

It is not hard to describe an affine subspace of $\R^n$ for any marked tubing which will contain the face of $\JG$ corresponding to that tubing.

\begin{defn} \label{d:nplanes}
Let $G$ be a graph with $n$ nodes  $v_i,\dots v_n$, and let $U$ be any marked tubing of $G$.  Let $\sm(v_i)$ be the smallest tube containing node $v_i.$  We define an affine subspace $H_U \subset \R^n$ by the following equations:
\begin{enumerate}
\item One equation for each thin tube $u$  given by:
$$\sum_{\sm(v_i)=u}x_i \ = \ 3^{|u|-1}-\sum_{s \Subset u }3^{|s|-1}.$$
\item One equation for each thick tube $u$  given by:
$$\sum_{\sm(v_i)=u}x_i \ = \ 3^{|u|}- \sum_{s \Subset u }3^{|s|}.$$
\end{enumerate}
\end{defn}

In the case of an upper or lower marked tubing $V$, the associated subspace $H_V$ is actually a \emph{hyperplane}, described by the single equation indicated by its single unbroken tube.  One result of Theorem~\ref{t:main} is that these are precisely the facet hyperplanes of a realization of the polytope $\JG$.  Indeed, the upper tubings correspond to facet-including hyperplanes that bound the polytope above, while the lower tubings yield facet-including hyperplanes that bound the polytope below.

\begin{exmp}
Figure~\ref{f:upperlower} shows examples of upper and lower tubings.  Based on the previous definition, the following hyperplanes will be associated to each of the appropriate tubings of the figure: \\ \\
\begin{tabular}{llll}
(a) \ $x_1 + x_2 + x_3 = 27$ \ \ & (b) \ $x_2 = 21$ & (c) \ $x_2 + x_3 = 24$ & \ \ \ (d) \ $x_1 = 21$ \\ \\
(e) \ $x_1 + x_2 + x_3 = 9$ &  (f) \ $x_1 + x_2 = 3$ \ \ & (g) \ $x_1 = 1$ \\ \\
\end{tabular}
\end{exmp}

\subsection{}

The multiplihedron  contains within its face structure several other important polytopes. The classic multiplihedron discovered by Stasheff, here corresponding to the graph multiplihedron of a path, encapsulates the combinatorics of a homotopy homomorphism between homotopy associative topological monoids. The important quotients of the Stasheff's multiplihedron then are the result of choosing a strictly associative domain or range for the maps to be studied. The case of a strictly associative range is described in \cite{sta2}, where Stasheff shows that the multiplihedron $\J_n$ becomes the associahedron $\K_{n+1}$.  The case of an associative domain is described in \cite{for2}, where the new quotient of the $n^{th}$ multiplihedron is the \emph{composihedron}, denoted $\mathcal{CK}(n)$.  These latter polytopes are the shapes of the axioms governing composition in higher enriched category theory, and thus referred to collectively as the composihedra. Finally the case of associativity of both range and domain is discussed in \cite{bv}, where the result is shown to be the $n$-dimensional cube.

Note that in Stasheff's multiplihedron, an associative domain corresponds to identifying certain points  within the lower facets, while an associative range corresponds to identifying certain points within the upper facets. In the case of a  graph multiplihedron, the simplest generalizations along these lines give rise to two families of convex polytopes.  We begin by demonstrating these two polytopes as convex hulls, using variations on Eq.~\eqref{e:real} which reflect the desired identifications.

\begin{defn}
The polytope $\JGd$ is the convex hull of the points $c_d(M_G)$ in $\R^n$ where
$$c_d:M_G \to \R^n: U \mapsto (x_1,\dots,x_n) \hspace{.5in} \textup{where} \ x_i =
\begin{cases}
\hspace{4.5mm} 1 & \text{if $\sm(v_i)$ is thin.} \\
3 \cdot f_U(v_i)  & \text{if $\sm(v_i)$ is thick.}
\end{cases} $$
This generalizes \emph{strict associativity of the domain} to graphs.
\end{defn}

A lower facet of $\JG$ is an isomorphic image of $\J \rec{t} \times \KG(t)$ for some thin tube $t.$  The quotient polytope $\JGd$ is achieved by identifying the images of any two points $(a,b) \sim (a,c)$ in  a lower facet, where $a$ is a point of $\J \rec{t}$ and $b,c$ are points in $\KG(t)$.  In terms of tubings,  the face poset of $\JGd$ is isomorphic to the poset $\JG$ modulo the equivalence relation on marked tubings generated by identifying any two tubings $U\sim V$ such that $U\prec V$ in $\JG$ precisely by the addition of a thin tube inside another thin tube, as in Figure 5(c).

\begin{defn}
The polytope $\JGr$ is the convex hull of the points $c_r(M_G)$ in $\R^n$ where
$$c_r:M_G \to \R^n: U \mapsto (x_1,\dots,x_n) \hspace{.5in} \textup{where} \ x_i =
\begin{cases}
 f_U(v_i) & \text{if $\sm(v_i)$ is thin.} \\
\hspace{4mm} 3^n  & \text{if $\sm(v_i)$ is thick.}
\end{cases} $$
This generalizes \emph{strict associativity of the range} to graphs.
\end{defn}

Recall that an upper facet of $\JG$ is the isomorphic image of
$$\K \rec{t_1\cup\dots\cup t_k} \times \JG(t_1) \times\dots\times\JG(t_k)$$
for broken tubes tube $t_1,\dots,t_k.$  The quotient polytope $\JGr$ is achieved by identifying the images of any two points $(x,y_1,\dots,y_k)\sim (z,y_1,\dots,y_k)$ in an upper facet.  In terms of tubings, the face poset of $\JGr$ is isomorphic to the poset $\JG$ modulo the equivalence relation on marked tubings generated by identifying any two tubings $U\sim V$ such that $U\prec V$ in $\JG$ precisely by the addition of a thick tube, as in Figure 5(e).

\begin{rem}
It is interesting to note that the polytope $\JGr$ appears in the context of deformations of bordered Riemann surfaces in \cite{dhv}, arising from the work of C.\ Liu \cite{liu}.  Indeed, it is the first example we know of where the associahedra (in the case when $G$ is a path) appear as truncations of cubes.
\end{rem}

Performing both quotienting operations simultaneously on the polytope $\JG$ will always yield the $n$-dimensional cube, where $n$ is the number of nodes of $G$.  Thus we have the following equation relating numbers of facets of the three polytopes defined here:
$$|\text{ facets of } \JGd | \ + \ |\text{ facets of } \JGr| \ - \ |\text{ facets of } \JG | \ = \ 2n$$
since the number of facets of the hypercube is $2n$.  Compare this with Lemma~\ref{l:enum}.  The following is a corollary of Theorem~\ref{t:perm} and the definitions above.  We leave it to the reader to fill in the details.

\begin{cor}
When $G$ is the complete graph, $\JGd$ and $\JGr$ are combinatorially equivalent.
\end{cor}

\begin{figure}[h]
\includegraphics[width=\textwidth]{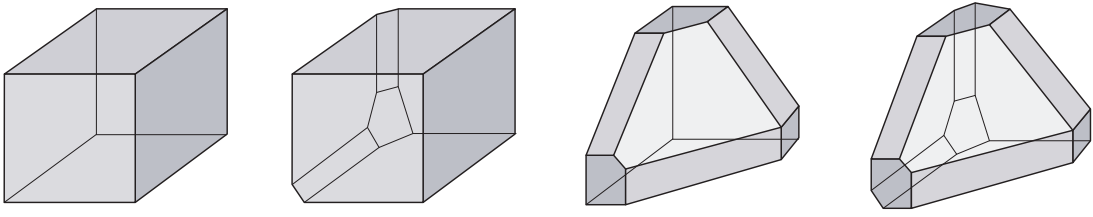}
\caption{A cube, associahedron, composihedron, and multiplihedron from a path with three nodes.}
\label{f:fromcubes}
\end{figure}

\begin{exmp}
When $G$ is a path with $n$ nodes, the polytopes $\JGd$ and $\JGr$ are the $n^{th}$ composihedron and the $(n+1)^{st}$ associahedron, respectively.  Figure~\ref{f:fromcubes} shows the realization of these polytopes discussed above for a path with three nodes.  Part (a) shows the cube, encapsulating associativity in both domain and range.  Parts (b) and (c) produce the associahedron and composihedron respectively; compare with \cite{dhv} and \cite{for2}.  Finally, truncating using the full collection of hyperplanes given by Eq.~\eqref{e:real} produces the multiplihedron.
\end{exmp}

%
%
\section{Proof of Theorem} \label{s:proof}
\subsection{}

The proof of Theorem~\ref{t:main} will use induction on the number of nodes of $G$.  This is feasible since we can characterize the structure
of the facets of our polytope via Propositions~\ref{p:lower} and~\ref{p:upper}.   Indeed, the dimension of the convex hull will be established, and together with the discovery of bounding hyperplanes, a characterization of the facets of the convex hull will be demonstrated.  Simultaneously, we can build a poset isomorphism out of (inductively assumed) isomorphisms that are restricted to the facets.  To begin, we will need a more general set of points and hyperplanes based on \emph{weights}.

Let $G$ be a graph with $n$ nodes, numbered by $i=1\dots n$.  Let $w_1,\dots,w_n$ be a list of positive integers (weights) which are associated to the respective nodes of $G$. For any tube $t$ of $G$, let
$$w(t) = \sum_{v_i\in t}w_i.$$
As before, let $\JG$ be the collection of maximal marked tubings of $G$.  Mimicking Eq.~\eqref{e:real}, we define a map from $M_G$ to $\R^n$ based on these weights.   Let $U$ be an element of $M_G$, $v$ a node of $G$, and $\sm(v)$ the smallest tube in $U$ containing $v$.  Let
$$f^w_U(v) = 3^{w(\sm(v))-1}-\sum_{s \Subset \sm(v) }3^{w(s)-1}$$
Now define $c^w : M_G \to \R^n$ where
\begin{equation}
c^w(U) = (x_1,\dots,x_n) \text{ where } x_i =
\begin{cases}
    \hspace{4.5mm} f^w_U(v_i) & \text{if  $\sm(v_i)$ thin} \\
    3 \cdot f^w_U(v_i) & \text{if $\sm(v_i)$ thick.}
\end{cases}
\label{e:weights}
\end{equation}

\begin{defn} \label{d:wnplanes}
Let $U$ be any marked tubing of $G$ and let $\sm(v_i)$ be the smallest tube containing node $v_i$.  Define an affine subspace $H^w_U \subset
\mathbb{R}^n$ by the following equations:
\begin{enumerate}
\item One equation for each  thin tube $u$ given by:
$$ \displaystyle{\sum_{\sm(v_i)=u}x_i = 3^{w(u)-1}-\sum_{s \Subset u }3^{w(s)-1}}$$
\item One equation for each thick tube $u$ given by:
$$\displaystyle{\sum_{\sm(v_i)=u}x_i = 3^{w(u)}- \sum_{s \Subset u }3^{w(s)}}.$$
\end{enumerate}
\end{defn}

\begin{lem} \label{l:w_graph_assoc}
Let $G$ be a graph with $n$ nodes.  Let $M^*_G$ be the subset of $M_G$ corresponding to all thick (or all thin) tubes.   The convex hull of $c^w(M^*_G)$ yields the graph associahedron $\KG$.
\end{lem}

\begin{proof}
This is seen by the remarks in Section 5 of \cite{dev2}. Having assigned weights $w_i$ to the nodes of $G,$ the function $\phi$ from unmarked tubes of $G$ to the integers is given by $\phi(u) = 3^{w(u)-1}$.  This function satisfies the inequality
$$\phi(u) > \phi(u_1) + \phi(u_2)$$
for $u_1, u_2$ any two proper subsets of the tube $u$.
To see this, let $W_1,\dots,W_k$ be the same list of weights as the $w_i$ for node $v_i$ in $u$, but ordered by decreasing size.
Let
$$W(u_j) = \sum_{i=1}^{|u_j|}W_i.$$
Now without loss of generality, let $|u_1| \ge |u_2|$.   Thus,
$$3^{w(u_1)-1} + 3^{w(u_2)-1} \leq 3^{W(u_1)-1} + 3^{W(u_1)-1} < 3 \cdot 3^{W(u_1)-1} \leq 3^{W(u)-1} =3^{w(u) -1},$$
as desired.
\end{proof}

\begin{prop} \label{p:dim}
For graph $G$ with $n$ nodes, the dimension of the convex hull of $c^w(M_G)$ is $n$.
\end{prop}

\begin{proof}
This is by inclusion of an $n$-dimensional prism within our convex hull.  For graph $G$ with $n$ nodes, there are two special tubings, the lower and  upper tubings for which the only tube is $G$.  Both of these, by Propositions~\ref{p:lower} and~\ref{p:upper}, have poset of refinements isomorphic to the unmarked tubings of $G$.  Thus, by Lemma~\ref{l:w_graph_assoc}, the convex hulls of the points associated to their respective maximal marked tubings are both isomorphic to the graph associahedra $\KG$ of dimension $n-1$.  Indeed, Eq.~\eqref{e:weights} shows the thick version scaled by a factor of three.  Moreover, the hyperplanes associated to these two tubings are parallel, so that the convex hull of just their vertices is a prism on $\KG$.  Thus, the entire dimension of $\JG$ is $n$.
\end{proof}

\subsection{}

The following three lemmas are needed for proving the main theorem.

\begin{lem} \label{l:onface}
Let $V$ be a facet of $\JG$ and let $U$ be a vertex
of $\JG$. If $U$ is a vertex of $V$, then $c^w(U)$ lies on $H^w_V$.
\end{lem}

\begin{proof}
First we note that if $V$ is a face of $U$, then $H^w_V \subset H^w_U$. This is true since when $V\prec U$ it implies that $V$ is obtained from $U$ by a sequence of any of possible moves described in Definition~\ref{d:poset}.  It is easily checked that each of these moves leaves inviolate the set of equations governing the coordinates $x_1,\dots,x_n$ induced by the original tubing, and introduces one new equation. The former is due to the fact that none of the refinements subtracts from the existing set of thick or thin tubes. The latter is due to the fact that each adds one more to the set of thick or thin tubes. Finally we point out that if $U$ is in $M_G$ then $H^w_U = c^w(U),$ since if $U$ is in $M_G$, then the tube
$u$ in each equation of Definition~\ref{d:nplanes} is the smallest tube containing $v_i$ for some node $v_i$.
\end{proof}

\begin{lem} \label{l:lowerface}
Let $V$ be a facet of $\JG$ which is a lower tubing and let $U$ be a vertex of $\JG$ such that $U \nprec V$.  Then $c^w(U)$ lies inside the halfspace of $\R^n$ created by $H^w_V$ not containing the origin.
\end{lem}

\begin{proof}
Let $t$ be the single thin tube of $V$.  For convenience, number the nodes of $G$ so $v_1, \ldots, v_k$ are the nodes of $t$ and let $c^w(U) = (x_1, \dots, x_n)$.  We must show for a vertex $U$ not in $V$,
$$x_1 + \dots + x_k \ > \ 3^{w(t) -1}.$$
This is seen by recognizing that $U$ either
\begin{enumerate}
\item contains a tube that is not compatible (as an unmarked tubing) with $t$, or
\item $\sm(v_i)$ is thick for some nodes $v_i$ of $t$.
\end{enumerate}
In the first case, there exists a node $v_i$ of $t$ for which $|\sm(v_i)| > |t|,$ and so $w(\sm(v_i)) > w(t).$  This leads to the desired inequality, regardless of the marking of $\sm(v_i)$.  If $t$ is compatible (as an unmarked tubing) with $U,$ but $\sm(v_i)$ is thick for some nodes $v_i$ of $t$,  the inequality follows simply due to the fact that $3 > 1$ (recall an additional factor of 3 in the definition of $c^w(U)$ for thick tubes).
\end{proof}

\begin{lem} \label{l:upperface}
Let $V$ be a facet of $\JG$ which is an upper tubing and let $U$  be a vertex of $\JG$ such that $U \nprec V$.  Then $c^w(U)$ lies inside the halfspace of $\R^n$ created by $H^w_V$ containing the origin.
\end{lem}

\begin{proof}
Let $t_1, \dots, t_r$ be the broken tubes of $V$.  For convenience, number the nodes of $G$ so $v_1,\dots,v_k$ are the nodes such that $\sm(v_i)=G$.  Let $c^w(U) = (x_1,\dots,x_n).$   We need to show for vertex $U$ not in $V$,
$$x_1 + \dots + x_k \ < \ 3^{w(G)}- \sum_{j=1}^r3^{w(t_j)}.$$
This is seen by recognizing that $U$ either
\begin{enumerate}
\item contains a tube that is not compatible (as an unmarked tubing) with some $t_j$, or
\item $\sm(v_i)$ is thin for some nodes $v_i$.
\end{enumerate}
For the first case, when all the tubes of $U$ are thick, the sum of all the coordinates in the point $c^w(U)$ is equal to $3^{w(G)}$. For some of the broken tubes $t_j$, there exists an included node $v_i$ for which $w(\sm(v_i)) > w(t_j).$  For these broken tubes, the sum of the coordinates calculated from nodes within $t_j$ is $3^{w(\sm(v_i))}.$  Thus,
$$x_1 + \dots + x_k \ = \ 3^{w(G)} - \sum_{j=1}^r \left(\sum_{v_i\in t_j}x_i\right) \ < \ 3^{w(G)} - \sum_{j=1}^r3^{w(t_j)},$$
since smaller terms are being subtracted in the last expression.  Again  if the underlying tubing of $U$ is preserved and more of the
tubes are allowed to be thin, the inequality is only strengthened. If $t_j$ is compatible (as an unmarked tubing) with $U$ but some nodes $v_i$  are such that  $\sm(v_i)$ is thin, then the inequality follows since $3 > 1.$
\end{proof}

\subsection{}

We are now in position to finish the proof of our key result.

\begin{proof}[Proof of Theorem~\ref{t:main}]
We will use induction on the number of nodes of $G$, made possible due to Propositions~\ref{p:lower} and \ref{p:upper}.
  We will proceed to prove that the theorem holds for the weighted version, with points $c^w(M_G)$ and that
  will imply the original version for all weights equal to 1.  The base case is when $n=1$.
  The two points in $\R^1$ are $3^{w_1-1}$ and $3^{w_1},$ whose convex hull is a line segment as expected.

The induction assumption is as follows:  For all graphs $G$ with number of nodes $k < n$  and for an arbitrary set of positive integer weights $w_1,\dots,w_k$, assume that the poset of marked tubings of $G$ is isomorphic to the face poset of the convex hull via the map $\varphi_G^w$ defined as follows:
$$\varphi_G^w: \JG \rightarrow CH\{c^w(M_G)\} : U' \mapsto CH\{c^w(U) ~|~ U \in M_G,~ U \prec U'\}.$$
Now we show this implies $\varphi_G^w$ to be an isomorphism in the case of $n$ nodes in $G$.

The mapping $\varphi_G^w$ clearly respects the ordering $\prec$ of marked tubings. This is evident since $U \prec U'$  implies for sets
$$\{V\in M_G ~|~ V \prec U\} \ \subset \ \{V'\in M_G ~|~ V' \prec U'\}.$$
Therefore the convex hulls obey the inclusion
$$CH\{c^w(V) ~ |~ V\in M_G ~,~ V \prec U\} \ \subset \ CH\{c^w(V') ~ |~ V'\in M_G ~,~ V' \prec U'\}.$$
Note that the restriction of $\varphi_G^w$ to tubings that are all thick (or all thin) is an isomorphism from the thick (thin) subposet to the face poset of the graph associahedra, by Lemma~\ref{l:w_graph_assoc}. We will denote these restrictions by $\varphi_G^w|_{thick}$ and $\varphi_G^w|_{thin}$ respectively.

Now by Propositions~\ref{p:lower} and \ref{p:upper}, the subposets of refinements of upper and lower tubings have the structure of cartesian products of tubing posets of certain smaller graphs.  This will allow the restriction of $\varphi^w_G$ to a lower or upper tubing $V$ to be shown to be an isomorphism:
$$\varphi_G^w|_V : \{U~|~U\prec V\} \to CH\{c^w(U)~|~ U \prec V,~ U \in M_G \}.$$
Keep in mind that the calculation of the coordinate $x_i$ is only affected by the structure of the tubing inside of the tube $\sm(v_i)$ which is the smallest tube containing node $v_i.$ Furthermore the calculation only reflects the size of the tubes $s \Subset \sm(v_i)$ and not their substructure.

For $V$ a lower tubing, with thin tube $t$, $\varphi_G^w|_V$ is an isomorphism since (up to renumbering of nodes)
$$\varphi_G^w|_V = \left(\varphi_{\rec{t}}^{\hat{w}} \times \varphi_{G(t)}^w|_{thin}\right)\circ\hat{\tau}^{-1},$$
where $\hat{\tau}$ is defined in Eq.~\eqref{e:tau}.
Each component of the first term is an isomorphism by induction. The new weighting $\hat{w}$ is determined by adding $w(t)$ to
each of the original weights
$w_i$ for which the node $v_i$ was connected to at least one node of $t$ by a single edge.

Similarly for upper tubes, the restriction of $\varphi_G$ to an upper tubing $V$ is given by (up to renumbering of nodes)
$$\varphi_G^w|_V = \left(\varphi_{\rec{t}}^{\hat{w}}|_{thick} \times \varphi_{G(t_1)}^w \times\dots\times\varphi_{G(t_k)}^w\right)\circ\hat{\eta}^{-1},$$
where $t$ is the union of broken tubes $t_1, \ldots, t_k$ and $\hat{\eta}$ is given by Eq.~\eqref{e:eta}. Each component of the first term is an isomorphism by induction. The new weight $\hat{w}$ is determined by adding $w(t_j)$ to each of the original weights $w_i$ for which the node $v_i$ was connected to at least one node of $t_j$ by a single edge.

\begin{nota}
We write $X \prec: Y$ if $X \prec Y$ and there does not exist a $Z$ such that $X \prec Z \prec Y$.
\end{nota}

\noindent $\bullet$ \  \emph{Show that $\varphi_G^w$ is injective}:

Let $X$ and $Y$ be two distinct marked tubings.  If $X,Y \prec V$, where $V \prec: G$, then  $\varphi_G^w(X) \ne  \varphi_G^w(Y)$ by induction, since as shown above the restriction of $\varphi_G^w$ to an upper or lower tubing is an isomorphism.  However, if $X \prec V\prec: G$ and $Y \nprec V$, then there exists $U \in M_G$ where $U \preccurlyeq Y$ and $U \nprec V.$  Then by Lemmas~\ref{l:lowerface} and \ref{l:upperface}, we have that $\varphi_G^w(X) \ne  \varphi_G^w(Y)$, since $c^w(U) \notin H_V^w$ and $\varphi_G^w(X) \subset H_V^w.$ 

\bigskip

\noindent $\bullet$ \  \emph{Show that $\varphi_G^w$ is surjective}:

The facets need to be described. First, we will show that the bounding hyperplanes $H_V^w$ each actually contain a facet of the convex hull. Then we will check that every facet is contained in one of these hyperplanes. The dimension of the facets  is now crucial. Recall that the total
dimension of the entire convex hull is $n$ by Proposition~\ref{p:dim}.  Now the dimension of the convex hull of the points associated to any upper or lower tubing is $n-1$ due to the following argument:  Since the dimension of $\JG$ is $n$ and the dimension of $\KG$ is $n-1,$ then the restriction of $\varphi_G^w$ to a lower tubing $V$ with thin tubing $t$ has image with dimension $(~n-|t|~) + (~|t|-1~) \ = \ n-1$. The restriction of $\varphi_G^w$ to an upper tubing $V$ with broken tubes $t_1,\dots, t_k$ has image with dimension $(~n - (~|t_1|+\dots +|t_k|~)) -1  + (~|t_1|+ \dots + |t_k|~) \ = \ n-1$ as well.

By Lemmas~\ref{l:lowerface} and \ref{l:upperface}, the hyperplanes $H_V^w$ are bounding planes that do contain the convex hulls of the restriction of $\varphi_G^w$ to the respective lower and upper tubings; thus, the image of that restriction is indeed a facet of the convex hull.  We now show the images of $\varphi_G^w|_V$ for the upper and lower tubings $V$ constitute the entire set of facets. This is equivalent to arguing that every codimension two face (a facet of the image of $\varphi_G^w|_V$) is also contained as a facet in $\varphi_G^w|_{V'}$ for some other upper or lower tubing $V'.$   By induction, the marked tubings $U\prec: V$ are the preimages of these codimension two faces. For each $U\prec: V$, it follows from  Definition~\ref{d:poset} that there is exactly one other upper or lower tubing $V'$ with $U\prec: V'$. Thus each codimension two face of the convex hull is contained in precisely two  of our set of upper and lower facets, showing that there can be no additional facets.

Finally, we prove that for any face $F$ of the the convex hull, there exists a tubing $U$ such that $\varphi_G^w(U) = F.$ If $F$ is a facet, we have already shown that $F = \varphi_G^w(V)$ for the corresponding upper or lower tubing V. Otherwise, let $F$ be a  convex hull of a collection of maximal marked tubings $\{U'\},$ and  $F\subset H_V^w$ for some upper or lower tubing $V.$  Then since $U' \prec V$ for each $U'$, there is a  preimage of $F$ by induction: the preimage of $F$ under $\varphi_G^w|_V.$
\end{proof}

\newpage

%
%
\bibliographystyle{amsplain}

\end{document}